\documentclass[12pt]{article}
\input{tcilatex}
\input{xy}

\xyoption{all}
\usepackage{graphicx}
\usepackage[utf8]{inputenc}
\usepackage{amsfonts,a4,epsf}

\begin{document}
\author{Rui Pedro Carpentier\footnote{rcarpentier@sapo.pt}\\Roger Picken\footnote{rpicken@math.ist.utl.pt}
\\
{\small\it Departamento de Matem\'{a}tica} \\ {\small\it Centro de
An\'{a}lise Matem\'{a}tica, Geometria e Sistemas Din\^{a}micos}\\
{\small\it Instituto Superior T\'{e}cnico}\\ {\small\it Avenida
Rovisco Pais, 1049-001 Lisboa}\\ {\small\it Portugal}}
\title{Some properties of Bowlin and Brin's color graphs}

\date{23rd April, 2018}
\maketitle
\begin{abstract}
Bowlin and Brin defined the class of color graphs, whose vertices are triangulated polygons compatible with a fixed four-coloring of the polygon vertices. In this article it is proven that each color graph has a vertex-induced embedding in a hypercube, and an upper bound is given for the hypercube dimension. The color graphs for $n$-gons up to $n=8$ are listed and studied, in particular enabling a question by Bowlin and Brin concerning the diameter of color graphs to be answered. Finally it is shown that color graphs with a certain type of subgraph cannot be isometrically embedded in a hypercube of any dimension.

\end{abstract}

\section{Introduction}
\label{int}

\begin{figure}[htbp]
\centering
\includegraphics[width=13cm]{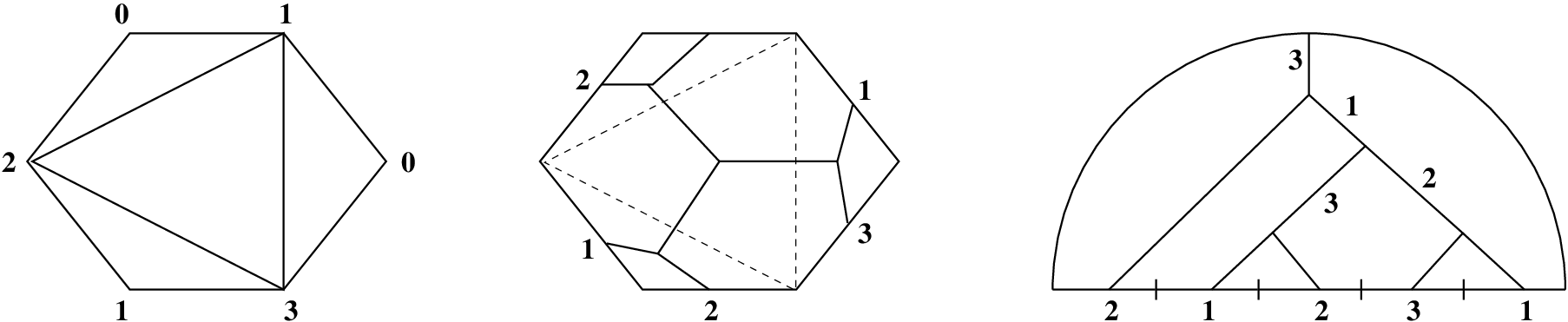}
\caption{From polygon triangulation to binary tree graph; polygon vertex coloring and the color vector ${\bf c}= 21231$ }
\label{fig:polygon-tree}
\end{figure}

\begin{figure}[htbp]
\centering
\includegraphics[width=13cm]{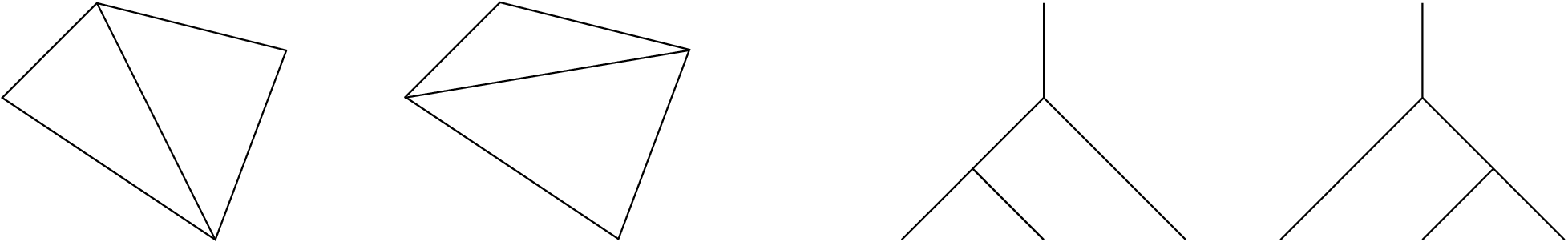}
\caption{Diagonal flip or rotation}
\label{fig:flip-rotation}
\end{figure}

In the context of their study of the four color theorem, Bowlin and Brin \cite{bb} introduced the notion of a color graph. We start by reviewing this definition.
Consider the set of triangulations of a regular polygon, such that the only vertices of each triangulation are those of the polygon itself. This set constitutes the vertex set of the associahedron of dimension $d$, denoted $A_d$, for a polygon with $d+3$ edges. Equivalently the vertices of $A_d$ are given by binary tree graphs with $d+2$ leaves. See Figure \ref{fig:polygon-tree} for the passage from a triangulation of a hexagon to the dual graph, and the deformation of the latter to display it as a tree graph with the root at the top and 5 leaves at the bottom (the numbering in the figure can be ignored at this stage). The edges of the associahedron are given by pairs of triangulations, or pairs of tree graphs, differing only by a single move, called a diagonal flip for triangulations or a rotation for tree graphs - see Figure \ref{fig:flip-rotation}. Both viewpoints have their advantages. Bowlin and Brin place more emphasis on the binary tree graph perspective, whereas we focus mainly on the triangulated polygon perspective, apart from in Theorem 1.

Consider the Klein four-group $V={\mathbb{Z}}_2\times {\mathbb{Z}}_2$, with elements taken to be $0,1,2,3$, such that $1+2+3=1+1=2+2=3+3=0$. A coloring of the vertices of a triangulated polygon with elements of $V$ is compatible with the triangulation if adjacent vertices are colored differently, as are vertices connected by a diagonal edge of the triangulation. See Figure \ref{fig:polygon-tree} for an example of a coloring of the vertices of a hexagon compatible with the triangulation shown. By assigning to each leaf of the dual tree graph the difference, or equivalently the sum,  of the colorings of the two vertices nearest to it, we obtain, reading from left to right, a list of non-zero elements of $V$, called a color vector by Bowlin and Brin \cite{bb}. Figure \ref{fig:polygon-tree} shows the color vector thus obtained in the example. The coloring of the leaves extends to a three-coloring of the edges of the whole tree graph, using only the colors 1, 2 and 3, by requiring that the sum of incident edge colors equals 0 at each internal vertex of the tree graph. Following Bowlin and Brin's terminology, such an edge coloring is termed proper, and we say that the color vector $\bf c$ is valid for the tree graph $T$ if it extends to a proper three-coloring of the edges of $T$, and that it is acceptable if it is valid for at least one tree graph $T$.

Given an acceptable color vector $\bf c$ consisting of $d+2$ elements of $V$ its color graph is the graph whose vertices are the vertices of $A_d$ such that $\bf c$ is valid for them, and whose edges are the edges of $A_d$ induced by these vertices. As pointed out by Bowlin and Brin, a result of \cite{gp} shows that any such color graph is either connected or has empty edge set (in the latter case the color vector is called rigid). The four color theorem is equivalent to the statement that for any pair of vertices of $A_d$ one can find at least one color graph to which they both belong. Triangulations in the same color graph can be connected by a sequence of signed diagonal flips, introduced by Eliahou \cite{e}, and further explored by Eliahou and Lecouvey \cite{el}. A criterium for determining whether a sequence of diagonal flips is a signed diagonal flip sequence was given by Carpentier \cite{c}. 

The main focus of Bowlin and Brin's work \cite{bb} was to explore the connection with Thompson's group $F$. However in the process they raised a number of questions concerning the properties of color graphs and gave partial answers. The purpose of the present article is to study further some properties of this interesting class of graphs, and to present a complete list of examples for color graphs associated to hexagons, heptagons and octagons.

\section{Embeddings of color graphs in integer lattices and hypercubes}

\begin{figure}[htbp]
\centering
\includegraphics[width=8cm]{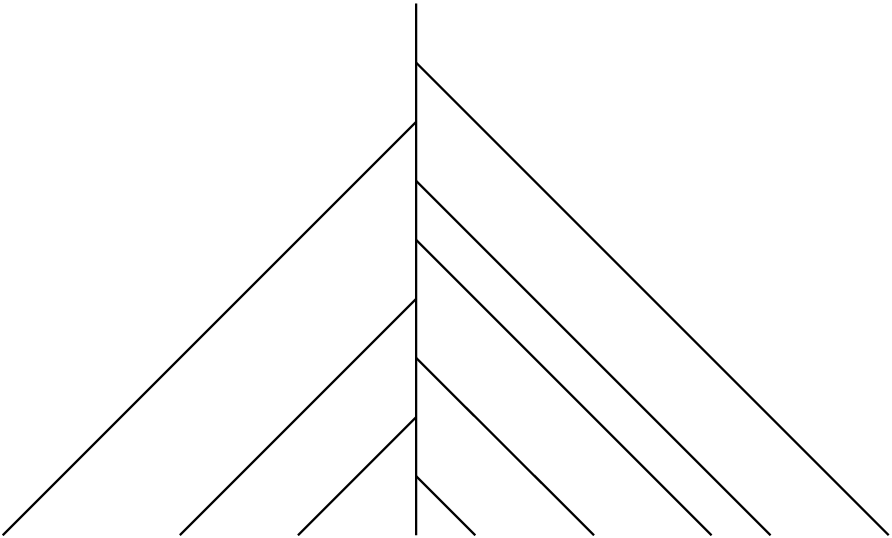}
\caption{Example of a vine associated to ${\bf c}=1^321^5$}
\label{fig:vine}
\end{figure}

The easiest color graphs to describe are those associated to color vectors of type $1^p21^q$, with $p,q\geq 0$. Any binary tree compatible with the color vector $1^p21^q$ is called a vine (see \cite{bb}), more properly a tree with a path from the root to leaf number $n+1$ (counting from the left) such that all the other leaves connect to that path (we will call this path the {\it central path}). See Figure \ref{fig:vine} for an example of a binary tree that is a vine. We will use the term {\it vine color graph} to describe color graphs associated to color vectors of the form $1^p21^q$, since every vertex of such a color graph is a vine. In \cite{bb} one of the questions asked was whether the color graphs with the greatest diameter within the class of $n$-gon color graphs, for a fixed value of $n$,  were always vine color graphs. See section \ref{exp} where we return to this question.

Let $L(m)$ be the integer lattice  graph of dimension $m$, i.e. the graph whose vertex set is the set of $m$-dimensional integer vectors and such that two vectors are adjacent iff they differ only in one entry and only by $\pm 1$.

\begin{theorem}
The color graph associated to the color vector $1^p21^q$ is isomorphic to the sub-graph of $L(m)$ induced by the vertex set $V=\{(x_1,\dots , x_m)\in \mathbb{Z}^m :0\leq x_1 \leq \cdots \leq x_m \leq q \}$.
\label{thm-vines}
\end{theorem}
\TeXButton{Proof}{\proof}
 There is a natural ordering of the leaves labelled with color 1: we say a leaf is above another leaf if it joins the central path above that other leaf. Label the leaves on the left from top to bottom by $i=1,\dots , m$, and let $x_i$ denote the number of leaves on the right that are above leaf $i$ on the left (in Figure \ref{fig:vine}, $x_1=1$, $x_2=3$, $x_3=4$).

We can then identify each such tree graph with a vector $(x_1,\dots ,x_m)\in L(m)$ with  $0\leq x_1 \leq \cdots \leq x_m \leq q$ and conversely, each such vector corresponds to a tree. Clearly two trees are connected by a rotation iff the corresponding vectors are adjacent in $L(m)$, since the order of precisely one leaf on the left and one on the right is altered by the rotation.
\TeXButton{End Proof}{\endproof}

The previous result leads us to ask in general when a color graph can be realized as a subgraph of $L(m)$ for some $m$, and what restrictions there are on the dimension $m$. It is well known that this is equivalent to asking whether a color graph can be realized as a subgraph of a hypercube of some finite dimension.

Given a finite set $S$ of cardinality $|S|=k$, the $k$-dimensional hypercube $\mathcal{H}(S)$ is defined to be a graph where the vertex set is $\mathcal{P}(S)$, the set of the subsets of $S$,  and two vertices $A$ and $B$ are connected by an edge of $\mathcal{H}(S)$ iff the symmetric difference $A\Delta B$ is a singleton.\footnote{The vertices of the hypercube can be identified with $\{0,1\}^S$ by associating to each subset its characteristic function, and under this identification adjacent vertices differ only in a single coordinate.} 

We start by showing that any color graph is a subgraph of a hypercube. 

\begin{theorem}
Every color graph is a vertex-induced subgraph of a hypercube.
\label{thm-hypercube}
\end{theorem}
\TeXButton{Proof}{\proof}
The color graph is determined by the vertex coloring of a polygon. A triangulation $T$ is compatible with the vertex coloring if all pairs of adjacent vertices on the edge of the polygon are colored differently, and likewise all pairs of vertices linked by a diagonal are colored differently. Choose one color to exclude (we call it $c$, the {\it forbidden color}), and consider the subset $S$ of diagonals connecting differently-colored vertices, neither of which is colored with $c$ (we call these the {\it allowed diagonals} and we call the diagonals that connect a vertex colored by $c$ with a vertex not colored by $c$ the {\it forbidden diagonals}).
Then for each triangulation $T$ of the polygon compatible with the coloring (i.e. for each vertex of the color graph) we assign the vertex of $\mathcal{H}(S)$, $D_T\subseteq S$, whose elements are  the allowed diagonals in the triangulation $T$. We have to show that this assignment is an injection on the vertex set of the color graph and that adjacent vertices are sent to adjacent vertices in the hypercube $\mathcal{H}(S)$. 

For the injectivity, we observe that we can recover the triangulation $T$ from its hypercube vertex $D_T$. Indeed, the allowed diagonals in $T$,  $D_T$, give a partial decomposition of the polygon into sub-polygons, and to recover the full triangulation $T$ we only need to identify the forbidden diagonals that were excluded from $T$. Every such sub-polygon with $k$ vertices, where $k> 3$, has exactly one vertex colored by $c$. 
It has to have at least one such vertex, since we need to use at least one forbidden diagonal to recover the original triangulation. If there were two or more vertices of the sub-polygon colored by $c$, these would have to be separated from each other in the original triangulation $T$, implying that there is at least one allowed diagonal inside the sub-polygon, which is a contradiction. Thus the forbidden diagonals in the sub-polygon are uniquely determined to be the $k-3$ diagonals which emanate from the single $c$-colored vertex.

To show that adjacent vertices, $T_1$ and $T_2$, are sent to adjacent vertices in the hypercube, we observe that in a flip precisely one diagonal of $T_1$ connecting vertices colored by two colors $c_1$ and $c_2$ is replaced by the complementary diagonal connecting vertices with the two complementary colors; one of these diagonals is forbidden and the other is allowed. Thus $D_{T_1}$ and  $D_{T_2}$ differ only in a single diagonal $\delta\in S$ and are therefore adjacent in $\mathcal{H}(S)$.

Finally, this subgraph of the hypercube is vertex-induced. If $D_{T_1}$ and  $D_{T_2}$ are adjacent in the hypercube differing only in a single diagonal $\delta\in S$ (say $D_{T_2}=D_{T_1}\setminus \{\delta\}$) then the diagonal $\delta$ separates a sub-polygon in the partial decomposition of $T_2$ into two sub-polygons in the partial decomposition of $T_1$, one with a vertex $v$ colored by $c$ and the other a triangle $t$ with no vertex colored by $c$. Thus $T_2$ is obtained from $T_1$ by performing the flip which exchanges the diagonal $\delta$ with the diagonal linking the vertex $v$ to the vertex of $t$ not incident to $\delta$. Thus the triangulations $T_1$ and $T_2$ are adjacent vertices in the color graph.
\TeXButton{End Proof}{\endproof}

As an immediate corollary we have:
\begin{corollary}
Every color graph is bipartite.
\label{cor:bipartite}
\end{corollary}

\vskip 0.2cm

In the proof of Theorem \ref{thm-vines}, we see that the dimension of the hypercube is given by the number of diagonals linking vertices of different colors neither of which is the forbidden color.

\begin{theorem}
A color graph associated to a polygon with $n$ sides can be embedded in a hypercube with dimension no greater than $\left\lfloor \frac{n(3n-8)}{16}\right\rfloor$.
\label{thm-hypercubedim}
\end{theorem}
\TeXButton{Proof}{\proof}
As in the proof of the previous theorem, we consider the vertex-coloring of a polygon defining the color graph, and then the dimension of the hypercube where the graph is embedded is given by the number of diagonals linking vertices of different colors neither of which is the forbidden color. Let $w$ be the number of vertices of the polygon colored by the forbidden color and $x$, $y$ and $z$ the number of vertices colored by the other three colors, thus $w+x+y+z=n$. 
Then the dimension of the hypercube would be the number of all diagonals (irrespective of the vertex colorings) minus the number of diagonals linking vertices of the same color minus the number of diagonals linking a vertex colored by $c$ to a vertex of a different color. Thus the dimension of the embedding is 
$$
\frac{n(n-3)}{2}-{w\choose 2}-{x\choose 2}-{y\choose 2}-{z\choose 2}-w(n-2-w) \quad\quad (\ast)
$$ 
Using
$$
\frac{n^2}{2} -{w\choose 2}-{x\choose 2}-{y\choose 2}-{z\choose 2} - w(n-w)= xy+yz+xz+ \frac{n}{2}
$$
the expression $(\ast)$ can be simplified to 
$$(x-1)(y-1)+(y-1)(z-1)+(z-1)(x-1)+n-3$$ 
This quadratic polynomial has a maximum value, subject to the constraints $x\ge 1$, $y\ge 1$, $z\ge 1$ and $x+y+z= n-w$, for $x=y=z=\frac{n-w}{3}$ and its maximum value is $\frac{(n-w-3)^2}{3}+n-3$.
Since we can always choose a forbidden color such that $w\ge \frac{n}{4}$, we get the upper bound $$\left\lfloor \frac{n(3n-8)}{16}\right\rfloor$$ for the dimension of the hypercube.
\TeXButton{End Proof}{\endproof}

\vskip 0.2cm
Another observation we can make is that a color graph only occupies a tiny portion of the hypercube. Since a triangulation of an $n$-gon has $n-3$ diagonals its representation on the hypercube has no more than $n-3$ non-null coordinates. On the other hand, choosing the forbidden color to be the one that colors most vertices, which is the case that minimizes the dimension of the hypercube, we have at least $\lceil \frac{n}{4}\rceil$ vertices with the forbidden color and in a triangulation they must be separated by $\lceil \frac{n}{4}\rceil -1$ allowed diagonals. Thus the immersion of the color graph in the hypercube lies in the region $\left\{A\subseteq S: \lceil \frac{n}{4}\rceil -1 \le |A| \le n-3\right\}$.

\section{Examples of color graphs}
\label{exp}

We have analysed all non-rigid color graphs corresponding to triangulations of hexagons, heptagons and octagons, and the detailed results are presented in the Appendix. In this section we explain our approach and point out some significant features of these examples.

\vskip 0.2cm

\noindent {\bf Listing of polygons and color graphs}

For each type of $n$-gon, we start by listing the possible four-colorings, up to permutations of colors and symmetries of the polygon. Each four-coloring implies a partition of $n$ into four integers, being the number of vertices of each color. We have grouped together the four-colorings with the same partition as follows:
\vskip 0.2cm

\noindent 5 hexagon colorings $P^6_1$ to $P^6_5$:\hskip 0.2cm 3,1,1,1 ($P^6_1$); 2,2,1,1 ($P^6_2$ to $P^6_5$)
\vskip 0.2cm

\noindent 7 heptagon colorings $P^7_1$ to $P^7_7$: \hskip 0.2cm 3,2,1,1 ($P^7_1$ to $P^7_3$); 2,2,2,1 ($P^7_4$ to $P^7_7$)
\vskip 0.2cm

\noindent 26 octagon colorings A-Z: 4,2,1,1 (A,B); 3,3,1,1 (C-G); 3,2,2,1 (H-S); 2,2,2,2 (T-Z)
\vskip 0.2cm

The color graphs themselves are listed as $G^6_1$ to $G^6_5$ (hexagon color graphs),  $G^7_1$ to $G^7_7$ (heptagon color graphs), and A-Z (octagon color graphs). 
\vskip 0.2cm

\noindent {\bf Embeddings in $L(m)$}

For the last four octagon color graphs W-Z we have also presented them as an embedding in $L(3)$, since we will observe (in section \ref{frc}) that all octagon color graphs can be embedded in $L(3)$ (for the other octagon color graphs this is obvious). In the same way all pentagon color graphs (just one graph) can be embedded in $L(1)$, and all hexagon and heptagon color graphs can be embedded in $L(2)$.
\vskip 0.2cm

\noindent {\bf Vines}

We identify the color graphs that are vine color graphs, i.e. the color graphs with a color vector of the form $1^p21^q$ where $p+q+2=n$ and we assume $p\leq q$ (see section \ref{int}). 
$$
\begin{array}{lc|lc|lc}
G^6_2 & 121^3 \quad & \quad G^7_1 & 121^4 \quad & \quad C & 121^5 \\
G^6_3 & 1^221^2 \quad & \quad G^7_4 & 1^221^3 \quad & \quad H & 1^221^4 \\
 & & & \quad & \quad T & 1^321^3
\end{array}
$$

\vskip 0.2cm

\noindent {\bf Diameter of color graphs}

In Bowlin-Brin \cite{bb}, Prop. 10.2, it was shown that the diameter of a vine color graph with color vector $1^p21^q$ is equal to $pq$, and they asked\footnote{We have reformulated and corrected the question in line with what was obviously the intention.} in Question 15.12 whether the vine color graphs with $p=q$ (for $n=2q+2$ even) or $p=q-1$  (for $n=2q+1$ odd) were always the color graphs with the greatest diameter for a given $n$. The examples show that, at least in some cases, the answer to this question is negative, since for heptagon color graphs ($n=7$) the diameter of the vine color graph $G^7_4$ is 6, which is exceeded by three other heptagon color graphs with diameters 7 or 8, and for octagon color graphs ($n=8$) the diameter of the vine color graph T is 9, which is exceeded by four other octagon color graphs, I, N, P, R, each with diameter 10. Despite these exceptions, the authors believe that the answer is affirmative for higher $n$.

\vskip 0.2cm
 
\noindent {\bf Size of lattice embeddings}

Apart from the dimension of the rectangular lattice in which the color graphs are embedded, we also make the observation that all octagon color graphs fit inside a rectangular region in $L(3)$ of dimensions $3\times 3\times 2$ (note the final figures for color graphs W,X,Y,Z in the appendix). Likewise heptagon and hexagon color graphs all fit inside a rectangular region in $L(2)$ of a fixed size, which can be $3\times 3$ or $4\times 2$ for heptagon color graphs, and $2\times 2$ or $3\times 1$ for hexagon color graphs.

\vskip 0.2cm

\noindent {\bf Inclusions and similarities of color graphs}

Apart from naturally finding many instances of color graphs for lower $n$ included as subgraphs of color graphs for higher $n$, we also note some cases of inclusions for the same $n$, e.g. for the three octagon color graphs T, U and X, we have that T is a subgraph of U and U is a subgraph of X. Likewise H is a subgraph of O, Q, V, W and Y. In other cases two color graphs for the same $n$ can be very similar in structure although neither is a subgraph of the other, e.g. the octagon color graphs P and S.

\vskip 0.2cm

\noindent {\bf Simplex shape of vine color graphs}

As a consequence of Theorem \ref{thm-vines} the vine color graphs associated to color vectors of the form $1^p2 1^q$ have the shape of simplices contained in $L(m)$. Thus, for instance, the color graphs $G^6_3$ and $G^7_4$ are shaped like triangles, and the octagon color graph T can be viewed as having a pyramid shape. We also note that vine color graphs with color vector of the form $121^q$ are always simple lines of length $q$, and can be embedded in $L(1)$.

\section{Further results and a couple of questions}
\label{frc}

In Theorem \ref{thm-hypercube} we showed that every color graph is a vertex-induced subgraph of a hypercube. We now ask whether this embedding is isometric. We note that this issue has also been studied from a general perspective by, amongst others, Djokovic \cite{d} and Ovchinnikov \cite{o}, although our discussion below is self-contained.

First of all we note that 
some octagon color graphs cannot be isometrically embedded in the hypercube of dimension 8 (being the dimension given in Theorem \ref{thm-hypercubedim}),  simply because their diameter is 9 or more. These are the octagon color graphs T, I, N, P, R.

Next we observe that some octagon color graphs (U, V, W, X, Y and Z) cannot be isometrically embedded in a hypercube of any dimension. This is a result of the following lemma and its corollary.

First we note, following \cite{o}, that the set $S$ underlying a hypercube $\mathcal{H}(S)$ may be infinite: the hypercube $\mathcal{H}(S)$ is then defined to be the graph where the vertex set is $\mathcal{P}_f(S)$, the set of all finite subsets of $S$, and as before, two vertices $A$ and $B$ are connected by an edge of $\mathcal{H}(S)$ iff the symmetric difference $A\Delta B$ is a singleton. The shortest path distance $d(A,B)$ on the hypercube $\mathcal{H}(S)$ is the {\it Hamming distance} between sets A and B: $d(A,B)=|A\Delta B|$ for $A,B\in \mathcal{P}_f(S)$.

Given two disjoint subsets  $T,C\subseteq S$ with $C$ finite, we define a face $\mathcal{F}_{T,C}$ on $\mathcal{H}(S)$ of dimension $k=|T|$ to be the subgraph induced by the vertex set $\{X\cup C: X\in \mathcal{P}_f(T)\}$.
Given a face $\mathcal{F}_{T,C}=\mathcal{F}$ on $\mathcal{H}(S)$ there is a natural projection map $\pi_\mathcal{F}$ from the vertices of $\mathcal{H}(S)$ to the vertices of $\mathcal{F}_{T,C}$ which sends a set $A$ to the set $A_{\mathcal{F}}:=(A\cap T)\cup C$. It is straightfoward to show, for any set $X\subseteq T$, that $A\Delta A_{\mathcal{F}}$ and $A_{\mathcal{F}}\Delta (X\cup C)$ are disjoint sets and $A\Delta (X\cup C)=(A\Delta A_{\mathcal{F}})\cup (A_{\mathcal{F}}\Delta (X\cup C))$. Thus we get the triangle equality $d(A,(X\cup C))=d(A,A_{\mathcal{F}})+d(A_{\mathcal{F}},(X\cup C))$, and therefore, $\pi_{\mathcal{F}}$ sends a vertex of $\mathcal{H}(S)$ to the unique vertex of $\mathcal{F}_{T,C}$ nearest to it.

\begin{lemma}
In a hypercube $\mathcal{H}(S)$, for any vertex $v$, we define the equivalence relation $\sim_v$ by saying that $x\sim_v y$ if $d(x,v)=d(y,v)$ ({\it i.e.} $x$ and $y$ are at the same distance from $v$). For a face $\mathcal{F}$ in $\mathcal{H}(S)$,  the equivalence relations  $\sim_v$ and $\sim_{\pi_{\mathcal{F}}(v)}$ are the same when restricted to the face $\mathcal{F}$. 
\end{lemma}
\TeXButton{Proof}{\proof}

The distance from $v$ to any vertex $x$ of $\mathcal{F}$ is equal to the sum of the distance from $v$ to $\pi_{\mathcal{F}}(v)$ and the distance from $\pi_{\mathcal{F}}(v)$ to $x$. Thus, any two vertices $x$ and $y$ of $F$ that are at the same distance from $v$ are at the same distance from $\pi_{\mathcal{F}}(v)$ and vice-versa, hence  the equivalence relations  $\sim_v$ and $\sim_{\pi_{\mathcal{F}}(v)}$ are the same when restricted to the face $\mathcal{F}$.
\TeXButton{End Proof}{\endproof}

\begin{corollary}
Let the {\it diamond ring graph} $DR_j$, $j\geq 1$, be the graph obtained from the 4-cycle graph by adding a path with length $j$ joining opposite vertices in the 4-cycle. A diamond ring graph $DR_{j}$ cannot be isometrically embedded in a hypercube of any dimension.
\end{corollary}
\TeXButton{Proof}{\proof}
If $j$ is odd this is obvious since then $DR_j$ is not bipartite. For $j=2k$ even, if $D_{2k}$ were isometrically embedded in a hypercube $C$, the 4-cycle of $D_{2k}$ would be a 2-dimensional face of $C$. Since the middle vertex $v$ of the length $2k$ path
 
is at the same distance ($k$ or $k+1$) from two opposite vertices in the 4-cycle, by the previous lemma there would be a vertex in the 4-cycle at the same distance from its opposite vertex as from itself, which is a contradiction.
\TeXButton{End Proof}{\endproof}
\vskip 0.2cm

We remark that this provides an alternative argument to that used by Ovchinnikov \cite{o}, Example 4.1, in a similar context.

Note that the above-mentioned octagon color graphs U, V, W, X, Y and Z, all have a diamond ring subgraph.

\vskip 0.2cm

We now look at the question of the minimum dimension for the integer lattice such that any $n$-gon color graph can be embedded in it. 
Looking at the examples in the appendix which were discussed in the previous section, we observe that any color graph associated to a coloring of an $n$-gon, for $n\leq 8$, is a subgraph (not necessarily vertex-induced) of the integer lattice graph of dimension $\lfloor \frac{n-2}{2}\rfloor$. It is easy to prove, as a consequence of theorem 1, that any vine color graph associated to a color vector of type $1^p21^q$ is a subgraph of the integer lattice graph of dimension $\min(p,q)$. So, for colorings of $n$-gons with these color vectors, their color graphs would be in $L(\lfloor \frac{n-2}{2}\rfloor)$. This raises the question of whether the {\it lattice dimension}\footnote{We define the lattice dimension for a color graph as the least dimension of a lattice graph where the color graph can be embedded. For a class of color graphs we define their lattice dimension as the maximum of the lattice dimensions of the graphs in the class.} of other color graphs follows the same pattern: $\lfloor \frac{n-2}{2}\rfloor$ for the class of $n$-gon color graphs. A couple of facts support this idea. First, the maximum degree of a color graph for an $n$-gon is no greater than $n-3$ which does not exceed the valency of the vertices in 
$L(\lfloor \frac{n-2}{2}\rfloor)$. And second, for an $n$-gon coloring its color graph contains an $r$-dimensional hypercube if there are $r$ non-overlapping four-colored quadrilaterals in the $n$-gon, therefore $r$ can not exceed $\lfloor \frac{n-2}{2}\rfloor$. However, it is possible that in higher dimensions the upper bound $\lfloor \frac{n-2}{2}\rfloor$ is not enough for the lattice dimension of a color graph. Because of theorem \ref{thm-hypercubedim} we know that the lattice dimension grows at most quadratically. We thus ask the following:

\vskip 0.2cm
\noindent {\bf Question 1} $\,$
Does the lattice dimension of $n$-gon color graphs grow linearly in $n$ or not?
\vskip 0.2cm

Our final question relates to Theorem \ref{thm-hypercubedim}:

\vskip 0.2cm
\noindent {\bf Question 2} $\,$ Is the dimension $\left\lfloor \frac{n(3n-8)}{16}\right\rfloor$ of Theorem \ref{thm-hypercubedim} the minimum dimension for 
which all $n$-gon color graphs have hypercube embeddings?
\vskip 0.2cm

\section*{Acknowledgements} This work was supported in part by the project UID/MAT/04459/2013 of the {\em Fundação para a Ciência e a Tecnologia} (FCT, Portugal).

\appendix
\section{Appendix.}

In this section we present a list of the color graphs associated to hexagons, heptagons and octagons. First, we list all possible four-colorings of a hexagon up to color permutations and hexagon symmetries:

$P^6_1$: \xymatrix@!0{
& *++[o][F-]{0} \ar@{-}[r]\ar@{-}[dl] & *++[o][F-]{1} \ar@{-}[dr] &  \\
*++[o][F-]{2}\ar@{-}[dr] & & & *++[o][F-]{0} \ar@{-}[dl] \\
& *++[o][F-]{0} \ar@{-}[r] & *++[o][F-]{3}  &  
}\quad $P^6_2$: \xymatrix@!0{
& *++[o][F-]{0} \ar@{-}[r]\ar@{-}[dl] & *++[o][F-]{1} \ar@{-}[dr] &  \\
*++[o][F-]{1}\ar@{-}[dr] & & & *++[o][F-]{0} \ar@{-}[dl] \\
& *++[o][F-]{2} \ar@{-}[r] & *++[o][F-]{3}  &  
}\quad $P^6_3$: \xymatrix@!0{
& *++[o][F-]{0} \ar@{-}[r]\ar@{-}[dl] & *++[o][F-]{1} \ar@{-}[dr] &  \\
*++[o][F-]{2}\ar@{-}[dr] & & & *++[o][F-]{0} \ar@{-}[dl] \\
& *++[o][F-]{3} \ar@{-}[r] & *++[o][F-]{2}  &  
}\quad 

$P^6_4$: \xymatrix@!0{
& *++[o][F-]{0} \ar@{-}[r]\ar@{-}[dl] & *++[o][F-]{1} \ar@{-}[dr] &  \\
*++[o][F-]{2}\ar@{-}[dr] & & & *++[o][F-]{0} \ar@{-}[dl] \\
& *++[o][F-]{1} \ar@{-}[r] & *++[o][F-]{3}  &  
}\quad $P^6_5$: \xymatrix@!0{
& *++[o][F-]{0} \ar@{-}[r]\ar@{-}[dl] & *++[o][F-]{1} \ar@{-}[dr] &  \\
*++[o][F-]{2}\ar@{-}[dr] & & & *++[o][F-]{3} \ar@{-}[dl] \\
& *++[o][F-]{1} \ar@{-}[r] & *++[o][F-]{0}  &  
}

with the corresponding color graphs:

\quad

$G^6_1$: \xymatrix{
\bullet \ar@{-}[d] &  \\
\bullet \ar@{-}[r]\ar@{-}[d] & \bullet \\
\bullet  &
} \quad $G^6_2$: \xymatrix{
\bullet \ar@{-}[d] &  \\
\bullet \ar@{-}[d] &  \\
\bullet \ar@{-}[r] & \bullet
} \quad $G^6_4$: \xymatrix{
\bullet \ar@{-}[d] & & \\
\bullet \ar@{-}[d] & & \\
\bullet \ar@{-}[r] & \bullet \ar@{-}[r] & \bullet
} $G^6_3$ and $G^6_5$: \xymatrix{
\bullet \ar@{-}[d] & & \\
\bullet \ar@{-}[d] \ar@{-}[r] & \bullet \ar@{-}[d] & \\
\bullet \ar@{-}[r] & \bullet \ar@{-}[r] & \bullet
} 

\quad

Next, we list all possible four-colorings of a heptagon up to color permutations and heptagon symmetries:

$P^7_1$: \xymatrix@!0{
& *++[o][F-]{0} \ar@{-}[r]\ar@{-}[dl] & *++[o][F-]{1} \ar@{-}[dr] &  \\
*++[o][F-]{1}\ar@{-}[d] & & & *++[o][F-]{0} \ar@{-}[d] \\
*++[o][F-]{0}\ar@{-}[dr] & & & *++[o][F-]{2} \ar@{-}[dll] \\
& *++[o][F-]{3}  &  &  
}\quad $P^7_2$: \xymatrix@!0{
& *++[o][F-]{0} \ar@{-}[r]\ar@{-}[dl] & *++[o][F-]{1} \ar@{-}[dr] &  \\
*++[o][F-]{2}\ar@{-}[d] & & & *++[o][F-]{0} \ar@{-}[d] \\
*++[o][F-]{0}\ar@{-}[dr] & & & *++[o][F-]{3} \ar@{-}[dll] \\
& *++[o][F-]{1}  &  &  
}\quad $P^7_3$: \xymatrix@!0{
& *++[o][F-]{0} \ar@{-}[r]\ar@{-}[dl] & *++[o][F-]{1} \ar@{-}[dr] &  \\
*++[o][F-]{2}\ar@{-}[d] & & & *++[o][F-]{0} \ar@{-}[d] \\
*++[o][F-]{0}\ar@{-}[dr] & & & *++[o][F-]{1} \ar@{-}[dll] \\
& *++[o][F-]{3}  &  &  
}

$P^7_4$: \xymatrix@!0{
& *++[o][F-]{0} \ar@{-}[r]\ar@{-}[dl] & *++[o][F-]{1} \ar@{-}[dr] &  \\
*++[o][F-]{1}\ar@{-}[d] & & & *++[o][F-]{0} \ar@{-}[d] \\
*++[o][F-]{2}\ar@{-}[dr] & & & *++[o][F-]{2} \ar@{-}[dll] \\
& *++[o][F-]{3}  &  &  
} \quad $P^7_5$: \xymatrix@!0{
& *++[o][F-]{0} \ar@{-}[r]\ar@{-}[dl] & *++[o][F-]{1} \ar@{-}[dr] &  \\
*++[o][F-]{2}\ar@{-}[d] & & & *++[o][F-]{0} \ar@{-}[d] \\
*++[o][F-]{1}\ar@{-}[dr] & & & *++[o][F-]{3} \ar@{-}[dll] \\
& *++[o][F-]{2}  &  &  
}

$P^7_6$: \xymatrix@!0{
& *++[o][F-]{0} \ar@{-}[r]\ar@{-}[dl] & *++[o][F-]{1} \ar@{-}[dr] &  \\
*++[o][F-]{2}\ar@{-}[d] & & & *++[o][F-]{0} \ar@{-}[d] \\
*++[o][F-]{1}\ar@{-}[dr] & & & *++[o][F-]{2} \ar@{-}[dll] \\
& *++[o][F-]{3}  &  &  
}\quad $P^7_7$: \xymatrix@!0{
& *++[o][F-]{0} \ar@{-}[r]\ar@{-}[dl] & *++[o][F-]{1} \ar@{-}[dr] &  \\
*++[o][F-]{2}\ar@{-}[d] & & & *++[o][F-]{2} \ar@{-}[d] \\
*++[o][F-]{1}\ar@{-}[dr] & & & *++[o][F-]{0} \ar@{-}[dll] \\
& *++[o][F-]{3}  &  &  
}

with the corresponding color graphs:

\quad

$G^7_1$: \xymatrix{
\bullet \ar@{-}[d] & &  \\
\bullet \ar@{-}[d] & &  \\
\bullet \ar@{-}[r] & \bullet \ar@{-}[r] & \bullet 
} \quad $G^7_2$: \xymatrix{
\bullet \ar@{-}[r] & \bullet \ar@{-}[d]\ar@{-}[r] & \bullet \ar@{-}[d] & \\
& \bullet \ar@{-}[d]\ar@{-}[r] & \bullet \ar@{-}[d]\ar@{-}[r] & \bullet \\
\bullet \ar@{-}[r] & \bullet \ar@{-}[r] & \bullet  &
} \quad $G^7_3$: \xymatrix{
& \bullet \ar@{-}[d] & & \bullet \ar@{-}[d] \\
\bullet \ar@{-}[r] & \bullet \ar@{-}[r] & \bullet \ar@{-}[r] & \bullet 
}   

$G^7_4$:\xymatrix{
\bullet \ar@{-}[d] & & & \\
\bullet \ar@{-}[r]\ar@{-}[d] &  \bullet \ar@{-}[d] & & \\
\bullet \ar@{-}[r]\ar@{-}[d] & \bullet \ar@{-}[r]\ar@{-}[d] &  \bullet \ar@{-}[d]  & \\
\bullet \ar@{-}[r] & \bullet \ar@{-}[r] & \bullet \ar@{-}[r] & \bullet   
} \quad $G^7_5$:\xymatrix{
\bullet \ar@{-}[d] & & & & \\
\bullet \ar@{-}[d] & & & & \\
\bullet \ar@{-}[d] & & & & \\
\bullet \ar@{-}[d] & & & & \\
\bullet \ar@{-}[r] & \bullet \ar@{-}[r] & \bullet \ar@{-}[r] & \bullet \ar@{-}[r] & \bullet 
}\quad

$G^7_6$:\xymatrix{
\bullet \ar@{-}[d] & & & & \\
\bullet \ar@{-}[r]\ar@{-}[d] & \bullet \ar@{-}[d] & & & \\
\bullet \ar@{-}[r]\ar@{-}[d] & \bullet \ar@{-}[r]\ar@{-}[d] &  \bullet \ar@{-}[d] & & \\
\bullet \ar@{-}[r] & \bullet \ar@{-}[r] & \bullet \ar@{-}[r] & \bullet \ar@{-}[r] & \bullet 
}\quad $G^7_7$:\xymatrix{
\bullet \ar@{-}[d] & & & & \\
\bullet \ar@{-}[d] & & & & \\
\bullet \ar@{-}[r]\ar@{-}[d] & \bullet \ar@{-}[d] & & & \\
\bullet \ar@{-}[r]\ar@{-}[d] & \bullet \ar@{-}[r]\ar@{-}[d] &  \bullet \ar@{-}[d] & & \\
\bullet \ar@{-}[r] & \bullet \ar@{-}[r] & \bullet \ar@{-}[r] & \bullet \ar@{-}[r] & \bullet 
}

\quad

And finally, we list all possible four-colorings of an octagon up to color permutations and octagon symmetries:

A: \xymatrix@!0{
& *++[o][F-]{0} \ar@{-}[r]\ar@{-}[dl] & *++[o][F-]{1} \ar@{-}[dr] &  \\
*++[o][F-]{2}\ar@{-}[d] & & & *++[o][F-]{0} \ar@{-}[d] \\
*++[o][F-]{0}\ar@{-}[dr] & & & *++[o][F-]{3} \ar@{-}[dl] \\
& *++[o][F-]{1} \ar@{-}[r] & *++[o][F-]{0}  &  
}\quad B: \xymatrix@!0{
& *++[o][F-]{0} \ar@{-}[r]\ar@{-}[dl] & *++[o][F-]{1} \ar@{-}[dr] &  \\
*++[o][F-]{1}\ar@{-}[d] & & & *++[o][F-]{0} \ar@{-}[d] \\
*++[o][F-]{0}\ar@{-}[dr] & & & *++[o][F-]{2} \ar@{-}[dl] \\
& *++[o][F-]{3} \ar@{-}[r] & *++[o][F-]{0}  &  
}\quad C: \xymatrix@!0{
& *++[o][F-]{0} \ar@{-}[r]\ar@{-}[dl] & *++[o][F-]{1} \ar@{-}[dr] &  \\
*++[o][F-]{1}\ar@{-}[d] & & & *++[o][F-]{0} \ar@{-}[d] \\
*++[o][F-]{0}\ar@{-}[dr] & & & *++[o][F-]{1} \ar@{-}[dl] \\
& *++[o][F-]{2} \ar@{-}[r] & *++[o][F-]{3}  &  
}\quad D: \xymatrix@!0{
& *++[o][F-]{0} \ar@{-}[r]\ar@{-}[dl] & *++[o][F-]{1} \ar@{-}[dr] &  \\
*++[o][F-]{1}\ar@{-}[d] & & & *++[o][F-]{0} \ar@{-}[d] \\
*++[o][F-]{0}\ar@{-}[dr] & & & *++[o][F-]{2} \ar@{-}[dl] \\
& *++[o][F-]{3} \ar@{-}[r] & *++[o][F-]{1}  &  
}

E: \xymatrix@!0{
& *++[o][F-]{0} \ar@{-}[r]\ar@{-}[dl] & *++[o][F-]{1} \ar@{-}[dr] &  \\
*++[o][F-]{2}\ar@{-}[d] & & & *++[o][F-]{0} \ar@{-}[d] \\
*++[o][F-]{0}\ar@{-}[dr] & & & *++[o][F-]{2} \ar@{-}[dl] \\
& *++[o][F-]{2} \ar@{-}[r] & *++[o][F-]{3}  &  
}\quad F: \xymatrix@!0{
& *++[o][F-]{0} \ar@{-}[r]\ar@{-}[dl] & *++[o][F-]{1} \ar@{-}[dr] &  \\
*++[o][F-]{2}\ar@{-}[d] & & & *++[o][F-]{2} \ar@{-}[d] \\
*++[o][F-]{0}\ar@{-}[dr] & & & *++[o][F-]{0} \ar@{-}[dl] \\
& *++[o][F-]{2} \ar@{-}[r] & *++[o][F-]{3}  &  
}\quad G: \xymatrix@!0{
& *++[o][F-]{0} \ar@{-}[r]\ar@{-}[dl] & *++[o][F-]{1} \ar@{-}[dr] &  \\
*++[o][F-]{2}\ar@{-}[d] & & & *++[o][F-]{2} \ar@{-}[d] \\
*++[o][F-]{0}\ar@{-}[dr] & & & *++[o][F-]{0} \ar@{-}[dl] \\
& *++[o][F-]{3} \ar@{-}[r] & *++[o][F-]{2}  &  
} \quad H: \xymatrix@!0{
& *++[o][F-]{0} \ar@{-}[r]\ar@{-}[dl] & *++[o][F-]{1} \ar@{-}[dr] &  \\
*++[o][F-]{1}\ar@{-}[d] & & & *++[o][F-]{0} \ar@{-}[d] \\
*++[o][F-]{0}\ar@{-}[dr] & & & *++[o][F-]{2} \ar@{-}[dl] \\
& *++[o][F-]{2} \ar@{-}[r] & *++[o][F-]{3}  &  
}

I: \xymatrix@!0{
& *++[o][F-]{0} \ar@{-}[r]\ar@{-}[dl] & *++[o][F-]{1} \ar@{-}[dr] &  \\
*++[o][F-]{2}\ar@{-}[d] & & & *++[o][F-]{0} \ar@{-}[d] \\
*++[o][F-]{0}\ar@{-}[dr] & & & *++[o][F-]{2} \ar@{-}[dl] \\
& *++[o][F-]{3} \ar@{-}[r] & *++[o][F-]{1}  &  
}\quad J: \xymatrix@!0{
& *++[o][F-]{0} \ar@{-}[r]\ar@{-}[dl] & *++[o][F-]{1} \ar@{-}[dr] &  \\
*++[o][F-]{2}\ar@{-}[d] & & & *++[o][F-]{0} \ar@{-}[d] \\
*++[o][F-]{0}\ar@{-}[dr] & & & *++[o][F-]{2} \ar@{-}[dl] \\
& *++[o][F-]{1} \ar@{-}[r] & *++[o][F-]{3}  &  
}\quad K: \xymatrix@!0{
& *++[o][F-]{0} \ar@{-}[r]\ar@{-}[dl] & *++[o][F-]{1} \ar@{-}[dr] &  \\
*++[o][F-]{2}\ar@{-}[d] & & & *++[o][F-]{0} \ar@{-}[d] \\
*++[o][F-]{0}\ar@{-}[dr] & & & *++[o][F-]{1} \ar@{-}[dl] \\
& *++[o][F-]{3} \ar@{-}[r] & *++[o][F-]{2}  &  
}\quad L: \xymatrix@!0{
& *++[o][F-]{0} \ar@{-}[r]\ar@{-}[dl] & *++[o][F-]{1} \ar@{-}[dr] &  \\
*++[o][F-]{2}\ar@{-}[d] & & & *++[o][F-]{0} \ar@{-}[d] \\
*++[o][F-]{0}\ar@{-}[dr] & & & *++[o][F-]{3} \ar@{-}[dl] \\
& *++[o][F-]{3} \ar@{-}[r] & *++[o][F-]{2}  &  
}

M: \xymatrix@!0{
& *++[o][F-]{0} \ar@{-}[r]\ar@{-}[dl] & *++[o][F-]{1} \ar@{-}[dr] &  \\
*++[o][F-]{2}\ar@{-}[d] & & & *++[o][F-]{0} \ar@{-}[d] \\
*++[o][F-]{0}\ar@{-}[dr] & & & *++[o][F-]{1} \ar@{-}[dl] \\
& *++[o][F-]{2} \ar@{-}[r] & *++[o][F-]{3}  &  
}\quad N: \xymatrix@!0{
& *++[o][F-]{0} \ar@{-}[r]\ar@{-}[dl] & *++[o][F-]{1} \ar@{-}[dr] &  \\
*++[o][F-]{1}\ar@{-}[d] & & & *++[o][F-]{2} \ar@{-}[d] \\
*++[o][F-]{0}\ar@{-}[dr] & & & *++[o][F-]{0} \ar@{-}[dl] \\
& *++[o][F-]{3} \ar@{-}[r] & *++[o][F-]{2}  &  
}\quad O: \xymatrix@!0{
& *++[o][F-]{0} \ar@{-}[r]\ar@{-}[dl] & *++[o][F-]{1} \ar@{-}[dr] &  \\
*++[o][F-]{1}\ar@{-}[d] & & & *++[o][F-]{2} \ar@{-}[d] \\
*++[o][F-]{0}\ar@{-}[dr] & & & *++[o][F-]{0} \ar@{-}[dl] \\
& *++[o][F-]{2} \ar@{-}[r] & *++[o][F-]{3}  &  
}\quad P: \xymatrix@!0{
& *++[o][F-]{0} \ar@{-}[r]\ar@{-}[dl] & *++[o][F-]{1} \ar@{-}[dr] &  \\
*++[o][F-]{2}\ar@{-}[d] & & & *++[o][F-]{2} \ar@{-}[d] \\
*++[o][F-]{0}\ar@{-}[dr] & & & *++[o][F-]{0} \ar@{-}[dl] \\
& *++[o][F-]{1} \ar@{-}[r] & *++[o][F-]{3}  &  
}

Q: \xymatrix@!0{
& *++[o][F-]{0} \ar@{-}[r]\ar@{-}[dl] & *++[o][F-]{1} \ar@{-}[dr] &  \\
*++[o][F-]{2}\ar@{-}[d] & & & *++[o][F-]{3} \ar@{-}[d] \\
*++[o][F-]{0}\ar@{-}[dr] & & & *++[o][F-]{0} \ar@{-}[dl] \\
& *++[o][F-]{1} \ar@{-}[r] & *++[o][F-]{3}  &  
}\quad R: \xymatrix@!0{
& *++[o][F-]{0} \ar@{-}[r]\ar@{-}[dl] & *++[o][F-]{1} \ar@{-}[dr] &  \\
*++[o][F-]{2}\ar@{-}[d] & & & *++[o][F-]{2} \ar@{-}[d] \\
*++[o][F-]{0}\ar@{-}[dr] & & & *++[o][F-]{0} \ar@{-}[dl] \\
& *++[o][F-]{3} \ar@{-}[r] & *++[o][F-]{1}  &  
}\quad S: \xymatrix@!0{
& *++[o][F-]{0} \ar@{-}[r]\ar@{-}[dl] & *++[o][F-]{1} \ar@{-}[dr] &  \\
*++[o][F-]{2}\ar@{-}[d] & & & *++[o][F-]{3} \ar@{-}[d] \\
*++[o][F-]{0}\ar@{-}[dr] & & & *++[o][F-]{0} \ar@{-}[dl] \\
& *++[o][F-]{3} \ar@{-}[r] & *++[o][F-]{1}  &  
}\quad T: \xymatrix@!0{
& *++[o][F-]{0} \ar@{-}[r]\ar@{-}[dl] & *++[o][F-]{1} \ar@{-}[dr] &  \\
*++[o][F-]{1}\ar@{-}[d] & & & *++[o][F-]{2} \ar@{-}[d] \\
*++[o][F-]{0}\ar@{-}[dr] & & & *++[o][F-]{3} \ar@{-}[dl] \\
& *++[o][F-]{3} \ar@{-}[r] & *++[o][F-]{2}  &  
}

U: \xymatrix@!0{
& *++[o][F-]{0} \ar@{-}[r]\ar@{-}[dl] & *++[o][F-]{1} \ar@{-}[dr] &  \\
*++[o][F-]{2}\ar@{-}[d] & & & *++[o][F-]{2} \ar@{-}[d] \\
*++[o][F-]{0}\ar@{-}[dr] & & & *++[o][F-]{3} \ar@{-}[dl] \\
& *++[o][F-]{3} \ar@{-}[r] & *++[o][F-]{1}  &  
}\quad V: \xymatrix@!0{
& *++[o][F-]{0} \ar@{-}[r]\ar@{-}[dl] & *++[o][F-]{1} \ar@{-}[dr] &  \\
*++[o][F-]{2}\ar@{-}[d] & & & *++[o][F-]{3} \ar@{-}[d] \\
*++[o][F-]{0}\ar@{-}[dr] & & & *++[o][F-]{2} \ar@{-}[dl] \\
& *++[o][F-]{1} \ar@{-}[r] & *++[o][F-]{3}  &  
}\quad W: \xymatrix@!0{
& *++[o][F-]{0} \ar@{-}[r]\ar@{-}[dl] & *++[o][F-]{1} \ar@{-}[dr] &  \\
*++[o][F-]{2}\ar@{-}[d] & & & *++[o][F-]{3} \ar@{-}[d] \\
*++[o][F-]{0}\ar@{-}[dr] & & & *++[o][F-]{2} \ar@{-}[dl] \\
& *++[o][F-]{3} \ar@{-}[r] & *++[o][F-]{1}  &  
}\quad X: \xymatrix@!0{
& *++[o][F-]{0} \ar@{-}[r]\ar@{-}[dl] & *++[o][F-]{1} \ar@{-}[dr] &  \\
*++[o][F-]{2}\ar@{-}[d] & & & *++[o][F-]{3} \ar@{-}[d] \\
*++[o][F-]{1}\ar@{-}[dr] & & & *++[o][F-]{0} \ar@{-}[dl] \\
& *++[o][F-]{3} \ar@{-}[r] & *++[o][F-]{2}  &  
}

Y: \xymatrix@!0{
& *++[o][F-]{0} \ar@{-}[r]\ar@{-}[dl] & *++[o][F-]{1} \ar@{-}[dr] &  \\
*++[o][F-]{2}\ar@{-}[d] & & & *++[o][F-]{3} \ar@{-}[d] \\
*++[o][F-]{1}\ar@{-}[dr] & & & *++[o][F-]{2} \ar@{-}[dl] \\
& *++[o][F-]{3} \ar@{-}[r] & *++[o][F-]{0}  &  
}\quad Z: \xymatrix@!0{
& *++[o][F-]{0} \ar@{-}[r]\ar@{-}[dl] & *++[o][F-]{1} \ar@{-}[dr] &  \\
*++[o][F-]{2}\ar@{-}[d] & & & *++[o][F-]{3} \ar@{-}[d] \\
*++[o][F-]{3}\ar@{-}[dr] & & & *++[o][F-]{2} \ar@{-}[dl] \\
& *++[o][F-]{1} \ar@{-}[r] & *++[o][F-]{0}  &  
}

\quad

with the corresponding color graphs:

\quad

A: \xymatrix{
& &\bullet \ar@{-}[d] & & \\
& \bullet \ar@{-}[r]\ar@{-}[d] & \bullet \ar@{-}[d]\ar@{-}[r] & \bullet \ar@{-}[d] & \\
\bullet \ar@{-}[r] & \bullet \ar@{-}[d]\ar@{-}[r] & \bullet \ar@{-}[d]\ar@{-}[r] & \bullet \ar@{-}[d]\ar@{-}[r] & \bullet \\
& \bullet \ar@{-}[r] & \bullet \ar@{-}[d]\ar@{-}[r] & \bullet  & \\
& & \bullet & &
} \quad B: \xymatrix{
\bullet \ar@{-}[d] & & \bullet \ar@{-}[d] \\
\bullet \ar@{-}[r]\ar@{-}[d] & \bullet \ar@{-}[r] & \bullet \ar@{-}[d] \\
\bullet  &  & \bullet 
} \quad C: \xymatrix{
\bullet \ar@{-}[d] & &  \\
\bullet \ar@{-}[d] & & \bullet \ar@{-}[d] \\
\bullet \ar@{-}[r] & \bullet \ar@{-}[r] & \bullet 
}

\quad

D: \xymatrix{
\bullet \ar@{-}[d] &  &  \\
\bullet \ar@{-}[d] & & \bullet \ar@{-}[d] \\
\bullet \ar@{-}[d] & & \bullet \ar@{-}[d] \\
\bullet \ar@{-}[r] & \bullet \ar@{-}[r]\ar@{-}[d] & \bullet \\
& \bullet & 
} \quad E and F: \xymatrix{
& \bullet \ar@{-}[d] \ar@{-}[r] & \bullet \ar@{-}[d] \ar@{-}[r] & \bullet \\
\bullet \ar@{-}[r] & \bullet \ar@{-}[d]\ar@{-}[r] & \bullet \ar@{-}[d] & \\
& \bullet \ar@{-}[d]\ar@{-}[r] & \bullet \ar@{-}[d]\ar@{-}[r] & \bullet \\
\bullet \ar@{-}[r] & \bullet \ar@{-}[r] & \bullet  &
}\quad G: \xymatrix@!0{
\bullet \ar@{-}[dd]\ar@{-}[rr]\ar@{-}[dr] && \bullet \ar@{-}[dd]\ar@{-}[rr] && \bullet \ar@{-}[dd]\ar@{-}[dr] &  \\
&\bullet &&&&\bullet \\
\bullet \ar@{-}[dd]\ar@{-}[rr] && \bullet \ar@{-}[dd]\ar@{-}[rr]\ar@{-}[dr] && \bullet \ar@{-}[dd] &  \\
&&& \bullet && \\
\bullet \ar@{-}[rr]\ar@{-}[dr]&& \bullet \ar@{-}[rr]&&\bullet \ar@{-}[dr]   \\
&\bullet &&&&\bullet 
}

\quad

H: \xymatrix{
\bullet \ar@{-}[d] & & & & \\
\bullet \ar@{-}[r]\ar@{-}[d] &  \bullet \ar@{-}[d] & & & \\
\bullet \ar@{-}[r]\ar@{-}[d] & \bullet \ar@{-}[r]\ar@{-}[d] &  \bullet \ar@{-}[d] & & \\
\bullet \ar@{-}[r]\ar@{-}[d] & \bullet \ar@{-}[r]\ar@{-}[d] & \bullet \ar@{-}[r]\ar@{-}[d] & \bullet \ar@{-}[d] & \\
\bullet \ar@{-}[r] & \bullet \ar@{-}[r] & \bullet \ar@{-}[r] & \bullet \ar@{-}[r] & \bullet 
} \quad I: \xymatrix@!0{
\bullet \ar@{-}[rr]\ar@{-}[dd] && \bullet \ar@{-}[dd]\ar@{-}[rr] && \bullet \ar@{-}[dd]\ar@{-}[rr] && \bullet \ar@{-}[dd]\ar@{-}[rr] && \bullet \ar@{-}[dd] \\
&&&&&&&&\\
\bullet  && \bullet \ar@{-}[dd]\ar@{-}[rr] && \bullet \ar@{-}[dd]\ar@{-}[rr]\ar@{-}[dl] && \bullet \ar@{-}[dd] && \bullet \ar@{-}[dd] \\
 &&& \bullet \ar@{-}[dd] &&&&& \\
 && \bullet \ar@{-}[r] &\ar@{-}[r]& \bullet \ar@{-}[rr]\ar@{-}[dl] && \bullet && \bullet \ar@{-}[dd]  \\
 &&& \bullet \ar@{-}[dd] &&&&& \\
&&&&&&&&\bullet\\
 &&& \bullet \ar@{-}[rr] && \bullet &&&  
}

\quad

J: \xymatrix@!0{
& \bullet \ar@{-}[dd]\ar@{-}[rr] && \bullet \ar@{-}[dd]\ar@{-}[rr]  && \bullet \ar@{-}[dd] &&\\
&&&&&&&\\
& \bullet  && \bullet \ar@{-}[d]\ar@{-}[rr]\ar@{-}[dl]  && \bullet \ar@{-}[dd]\ar@{-}[dl]\ar@{-}[rr] && \bullet \ar@{-}[dd] \\
\bullet \ar@{-}[dd]\ar@{-}[rr]  && \bullet \ar@{-}[dd]\ar@{-}[rr] &\ar@{-}[d]& \bullet \ar@{-}[dd] &&& \\
& && \bullet \ar@{-}[r]\ar@{-}[dl]  &\ar@{-}[r]& \bullet \ar@{-}[dl]\ar@{-}[rr] && \bullet \ar@{-}[dd] \\
\bullet \ar@{-}[dd]\ar@{-}[rr]  && \bullet \ar@{-}[dd]\ar@{-}[rr] && \bullet \ar@{-}[dd] &&& \\
& \bullet \ar@{-}[dl] &&&&\bullet \ar@{-}[rr] &&\bullet \\
\bullet \ar@{-}[rr]  && \bullet \ar@{-}[rr] && \bullet &&& 
} \quad K: \xymatrix@!0{
\bullet \ar@{-}[rr] && \bullet \ar@{-}[dd]\ar@{-}[rr] && \bullet \ar@{-}[dd]\ar@{-}[rr] && \bullet \ar@{-}[dd]\ar@{-}[rr] && \bullet  \\
&&&&&&&&\\
&& \bullet \ar@{-}[dd]\ar@{-}[rr] && \bullet \ar@{-}[dd]\ar@{-}[rr]\ar@{-}[dl] && \bullet \ar@{-}[dd] &&  \\
 &&& \bullet \ar@{-}[dd] &&&&& \\
 && \bullet \ar@{-}[r] &\ar@{-}[r]& \bullet \ar@{-}[rr]\ar@{-}[dl] && \bullet &&   \\
 &&& \bullet \ar@{-}[dd] &&&&& \\
&&&&&&&&\\
 &&& \bullet \ar@{-}[rr] && \bullet \ar@{-}[rr] && \bullet &  
}

\quad

L: \xymatrix{
\bullet \ar@{-}[d] & & & & & & &  \\
\bullet \ar@{-}[d]\ar@{-}[r] & \bullet \ar@{-}[r] & \bullet \ar@{-}[d]\ar@{-}[r] & \bullet \ar@{-}[d]\ar@{-}[r] & \bullet \ar@{-}[d]\ar@{-}[r] & \bullet \ar@{-}[d]\ar@{-}[r] & \bullet \ar@{-}[d]\ar@{-}[r] & \bullet  \\ 
\bullet  &  & \bullet \ar@{-}[r] & \bullet \ar@{-}[d]\ar@{-}[r] & \bullet \ar@{-}[d]\ar@{-}[r] & \bullet \ar@{-}[d]\ar@{-}[r] & \bullet  &  \\ 
  &  &  & \bullet \ar@{-}[r] & \bullet \ar@{-}[d]\ar@{-}[r] & \bullet  &  &  \\ 
  &  &  &  & \bullet &  &  &  
} \quad M: \xymatrix{
& \bullet \ar@{-}[d] & \\
\bullet \ar@{-}[d] \ar@{-}[r] & \bullet \ar@{-}[r] & \bullet \ar@{-}[d] \\
\bullet \ar@{-}[d] &  & \bullet \ar@{-}[d] \\
\bullet \ar@{-}[d] &  & \bullet \ar@{-}[d] \\
\bullet  &  & \bullet 
}

\quad

N: \xymatrix{
\bullet \ar@{-}[d]\ar@{-}[r] & \bullet \ar@{-}[r] & \bullet \ar@{-}[d]\ar@{-}[r] & \bullet \ar@{-}[r]& \bullet \ar@{-}[d] \\
\bullet &  & \bullet & & \bullet \ar@{-}[d] \\
 & & & & \bullet \ar@{-}[d] \\
 & & \bullet \ar@{-}[r] & \bullet \ar@{-}[r] & \bullet 
} \quad O: \xymatrix{
 & & & \bullet \ar@{-}[d] & & \\
\bullet \ar@{-}[d]\ar@{-}[r] & \bullet \ar@{-}[d]\ar@{-}[r] & \bullet \ar@{-}[d]\ar@{-}[r] & \bullet \ar@{-}[d]\ar@{-}[r] & \bullet \ar@{-}[d] & \\
 \bullet \ar@{-}[d]\ar@{-}[r] & \bullet \ar@{-}[d]\ar@{-}[r] & \bullet \ar@{-}[d]\ar@{-}[r] & \bullet \ar@{-}[r] & \bullet \ar@{-}[r] & \bullet \\
\bullet \ar@{-}[r]\ar@{-}[d] & \bullet \ar@{-}[d]\ar@{-}[r] & \bullet & & & \\
\bullet \ar@{-}[r]\ar@{-}[d] & \bullet & & & & \\
\bullet & & & & & 
}

\quad

P: \xymatrix@!0{
&  && \bullet \ar@{-}[dd]\ar@{-}[rr] && \bullet \ar@{-}[dd]\ar@{-}[rr] && \bullet \ar@{-}[dd]\ar@{-}[rr] && \bullet \ar@{-}[dd]\\
\bullet \ar@{-}[dd] &&&&&&&&&\\
& && \bullet \ar@{-}[d]\ar@{-}[rr]\ar@{-}[dl]  && \bullet \ar@{-}[dd]\ar@{-}[dl]\ar@{-}[rr] && \bullet \ar@{-}[dd] && \bullet\\
\bullet \ar@{-}[dd]\ar@{-}[rr]  && \bullet \ar@{-}[dd]\ar@{-}[rr] &\ar@{-}[d]& \bullet \ar@{-}[dd] &&&&& \\
& && \bullet \ar@{-}[r]\ar@{-}[dl]  &\ar@{-}[r]& \bullet \ar@{-}[dl]\ar@{-}[rr] && \bullet && \\
\bullet \ar@{-}[rr]  && \bullet \ar@{-}[dd]\ar@{-}[rr] && \bullet \ar@{-}[dd] &&&& \bullet \ar@{-}[dd] & \\
& &&&&&&&& \\
&& \bullet \ar@{-}[rr] && \bullet \ar@{-}[rr] && \bullet \ar@{-}[rr] && \bullet & 
}\quad Q: \xymatrix{
 & \bullet \ar@{-}[d] & & & & \\
\bullet \ar@{-}[d]\ar@{-}[r] & \bullet \ar@{-}[d]\ar@{-}[r] & \bullet \ar@{-}[d]\ar@{-}[r] & \bullet \ar@{-}[d]\ar@{-}[r] & \bullet \ar@{-}[d] & \\
\bullet \ar@{-}[r]\ar@{-}[d] & \bullet \ar@{-}[r] & \bullet \ar@{-}[d]\ar@{-}[r] & \bullet \ar@{-}[d]\ar@{-}[r] & \bullet \ar@{-}[d] & \\
\bullet &  & \bullet \ar@{-}[r] & \bullet \ar@{-}[d]\ar@{-}[r] & \bullet \ar@{-}[d] & \\
 &  &  & \bullet \ar@{-}[d]\ar@{-}[r] & \bullet \ar@{-}[r]\ar@{-}[d] & \bullet\\
 &  &  & \bullet \ar@{-}[d]\ar@{-}[r] & \bullet  & \\
 &  &  & \bullet  &  & \\
}

R: \xymatrix{
\bullet \ar@{-}[d] & & & \\ 
\bullet \ar@{-}[d]\ar@{-}[r] & \bullet \ar@{-}[d] & &  \\ 
\bullet \ar@{-}[d]\ar@{-}[r] & \bullet \ar@{-}[d]\ar@{-}[r] & \bullet \ar@{-}[d] &  \\
\bullet \ar@{-}[d]\ar@{-}[r] & \bullet \ar@{-}[d]\ar@{-}[r] & \bullet \ar@{-}[d]\ar@{-}[r] & \bullet \ar@{-}[d] \\
\bullet \ar@{-}[d]\ar@{-}[r] & \bullet \ar@{-}[d]\ar@{-}[r] & \bullet & \bullet \ar@{-}[d] \\
\bullet \ar@{-}[d]\ar@{-}[r] & \bullet  &  & \bullet  \\
\bullet \ar@{-}[r] & \bullet \ar@{-}[r] & \bullet &  
}\quad S: \xymatrix@!0{
&&&&&&&&& \bullet \ar@{-}[dd]\\
&&&&&&&&& \\
&  && \bullet \ar@{-}[dd]\ar@{-}[rr] && \bullet \ar@{-}[dd]\ar@{-}[rr] && \bullet \ar@{-}[dd]\ar@{-}[rr] && \bullet \ar@{-}[dd]\\
\bullet \ar@{-}[dd] &&&&&&&&&\\
& && \bullet \ar@{-}[d]\ar@{-}[rr]\ar@{-}[dl]  && \bullet \ar@{-}[dd]\ar@{-}[dl]\ar@{-}[rr] && \bullet \ar@{-}[dd] && \bullet\\
\bullet \ar@{-}[dd]\ar@{-}[rr]  && \bullet \ar@{-}[dd]\ar@{-}[rr] &\ar@{-}[d]& \bullet \ar@{-}[dd] &&&&& \\
& && \bullet \ar@{-}[r]\ar@{-}[dl]  &\ar@{-}[r]& \bullet \ar@{-}[dl]\ar@{-}[rr] && \bullet && \\
\bullet \ar@{-}[rr]  && \bullet \ar@{-}[dd]\ar@{-}[rr] && \bullet \ar@{-}[dd] &&&&  & \\
& &&&&&&&& \\
&& \bullet \ar@{-}[rr] && \bullet \ar@{-}[rr] && \bullet  && & 
}

T: \xymatrix@!0{
&  && \bullet \ar@{-}[dd]\ar@{-}[rr] && \bullet \ar@{-}[dd]\ar@{-}[rr] && \bullet \ar@{-}[dd]\ar@{-}[dl] \\
&&&&&&\bullet & \\
& && \bullet \ar@{-}[d]\ar@{-}[rr]\ar@{-}[dl]  && \bullet \ar@{-}[dd]\ar@{-}[dl]\ar@{-}[rr] && \bullet \ar@{-}[dd] \\
\bullet \ar@{-}[dd]\ar@{-}[rr]  && \bullet \ar@{-}[dd]\ar@{-}[rr] &\ar@{-}[d]& \bullet \ar@{-}[dd] &&& \\
& && \bullet \ar@{-}[r]\ar@{-}[dl]  &\ar@{-}[r]& \bullet \ar@{-}[dl]\ar@{-}[rr] && \bullet \\
\bullet \ar@{-}[rr]\ar@{-}[dd]  && \bullet \ar@{-}[dd]\ar@{-}[rr] && \bullet \ar@{-}[dd] &&& \\
& \bullet &&&&&& \\
\bullet \ar@{-}[ru]\ar@{-}[rr] && \bullet \ar@{-}[rr] && \bullet  && & 
}\quad U: \xymatrix@!0{
& \bullet \ar@{-}[dl]\ar@{-}[d] && \bullet \ar@{-}[d]\ar@{-}[rr] && \bullet \ar@{-}[d]\ar@{-}[rr] && \bullet \ar@{-}[dd]\ar@{-}[dl]  \\
\bullet \ar@{-}[rr] &\ar@{-}[d]& \bullet \ar@{-}[rr] &\ar@{-}[d]& \bullet \ar@{-}[rr] &\ar@{-}[d]&\bullet & \\
& \bullet \ar@{-}[d] && \bullet \ar@{-}[d]\ar@{-}[rr]\ar@{-}[dl]  && \bullet \ar@{-}[dd]\ar@{-}[dl]\ar@{-}[rr] && \bullet \ar@{-}[dd]   \\
\bullet \ar@{-}[dd]\ar@{-}[rr]  & \ar@{-}[d] & \bullet \ar@{-}[dd]\ar@{-}[rr] &\ar@{-}[d]& \bullet \ar@{-}[dd] &&& \\
& \bullet \ar@{-}[d] && \bullet \ar@{-}[r]\ar@{-}[dl]  &\ar@{-}[r]& \bullet \ar@{-}[dl]\ar@{-}[rr] && \bullet \\
\bullet \ar@{-}[rr]\ar@{-}[dd]  & \ar@{-}[d] & \bullet \ar@{-}[dd]\ar@{-}[rr] && \bullet \ar@{-}[dd] &&& \\
& \bullet &&&&&& \\
\bullet \ar@{-}[ru]\ar@{-}[rr] && \bullet \ar@{-}[rr] && \bullet  && &
}

\quad

V: \xymatrix@!0{
&&& \bullet \ar@{-}[dl]\ar@{-}[rr] && \bullet \ar@{-}[dl]\ar@{-}[rr] && \bullet \ar@{-}[dl]\ar@{-}[rr] && \bullet \ar@{-}[dl]\ar@{-}[dd]  \\
&& \bullet \ar@{-}[dl]\ar@{-}[rr] && \bullet \ar@{-}[dl]\ar@{-}[rr] && \bullet \ar@{-}[dl]\ar@{-}[rr] && \bullet  & \\
& \bullet \ar@{-}[dl]\ar@{-}[rr] && \bullet \ar@{-}[dl]\ar@{-}[rr] && \bullet  &&&& \bullet  \ar@{-}[dd] \\
\bullet \ar@{-}[dd]\ar@{-}[rr] && \bullet  &&&&&&& \\
&&& \bullet \ar@{-}[dl]\ar@{-}[rr] && \bullet \ar@{-}[dl]\ar@{-}[rr] && \bullet \ar@{-}[dl]\ar@{-}[rr] && \bullet \ar@{-}[dl]  \\
\bullet \ar@{-}[dd]&& \bullet \ar@{-}[dl]\ar@{-}[rr] && \bullet \ar@{-}[dl]\ar@{-}[rr] && \bullet \ar@{-}[dl]\ar@{-}[rr] && \bullet  & \\
& \bullet \ar@{-}[dl]\ar@{-}[rr] && \bullet \ar@{-}[dl]\ar@{-}[rr] && \bullet  &&&&  \\
\bullet \ar@{-}[rr] && \bullet  &&&&&&& 
}\quad

W: \xymatrix{
 & & & \bullet \ar@{-}[d]\ar@{-}[r] & \bullet \ar@{-}[d]\ar@{-}[r] & \bullet \ar@{-}[d]\ar@{-}[r] & \bullet \ar@{-}[d]\ar@{-}[r] & \bullet \ar@{-}[ddddddd] \\
 & & & \bullet \ar@{-}[d]\ar@{-}[r] & \bullet \ar@{-}[d]\ar@{-}[r] & \bullet \ar@{-}[d]\ar@{-}[r] & \bullet &  \\
 & & & \bullet \ar@{-}[d]\ar@{-}[r] & \bullet \ar@{-}[d]\ar@{-}[r] & \bullet  & & \\
\bullet \ar@{-}[d]\ar@{-}[r] & \bullet \ar@{-}[d]\ar@{-}[r] & \bullet \ar@{-}[d]\ar@{-}[r] & \bullet \ar@{-}[d]\ar@{-}[r] & \bullet \ar@{-}[d] &  & & \\
\bullet \ar@{-}[d]\ar@{-}[r] & \bullet \ar@{-}[d]\ar@{-}[r] & \bullet \ar@{-}[d]\ar@{-}[r] & \bullet \ar@{-}[r] & \bullet  &  & & \\
\bullet \ar@{-}[d]\ar@{-}[r] & \bullet \ar@{-}[d]\ar@{-}[r] & \bullet  & &  &  & & \\
\bullet \ar@{-}[d]\ar@{-}[r] & \bullet  & & &  &  & & \\
\bullet \ar@{-}[rrrrrrr] &  & & &  & & & \bullet \\
}

\quad

X: \xymatrix@!0{
&  && \bullet \ar@{-}[dd]\ar@{-}[rr] && \bullet \ar@{-}[dd]\ar@{-}[rr] && \bullet \ar@{-}[dd]\ar@{-}[drrrr] &&&&&\\
&&&&&&& &&&& \bullet \ar@{-}[dd] &\\
& && \bullet \ar@{-}[d]\ar@{-}[rr]\ar@{-}[dl]  && \bullet \ar@{-}[dd]\ar@{-}[dl]\ar@{-}[rr] && \bullet \ar@{-}[dd] & \bullet \ar@{-}[dd]\ar@{-}[rr] && \bullet \ar@{-}[dd]\ar@{-}[ru] && \\
 && \bullet \ar@{-}[dd]\ar@{-}[rr] &\ar@{-}[d]& \bullet \ar@{-}[dd] &&& &&&& \bullet \ar@{-}[dd] &\\
& \bullet \ar@{-}[dd]\ar@{-}[ru] && \bullet \ar@{-}[r]\ar@{-}[dl]  &\ar@{-}[r]& \bullet \ar@{-}[dl]\ar@{-}[rr] && \bullet & \bullet \ar@{-}[d]\ar@{-}[rr]\ar@{-}[dl]  && \bullet \ar@{-}[dd]\ar@{-}[dl]\ar@{-}[ru] && \\
 && \bullet \ar@{-}[dd]\ar@{-}[rr] && \bullet \ar@{-}[dd] & \bullet \ar@{-}[dd]\ar@{-}[rr]  && \bullet \ar@{-}[dd]\ar@{-}[rr] &\ar@{-}[d]& \bullet \ar@{-}[dd] && \bullet &  \\
& \bullet \ar@{-}[dd]\ar@{-}[ru] &&&&&& & \bullet \ar@{-}[r]\ar@{-}[dl]  &\ar@{-}[r]& \bullet \ar@{-}[dl]\ar@{-}[ru] && \\
&& \bullet \ar@{-}[rr] && \bullet  & \bullet \ar@{-}[rr]\ar@{-}[dd]  && \bullet \ar@{-}[dd]\ar@{-}[rr] && \bullet \ar@{-}[dd] &&& \\
& \bullet \ar@{-}[drrrr]\ar@{-}[ru] &&&&&& &&&&&\\
&&&&& \bullet \ar@{-}[rr] && \bullet \ar@{-}[rr] && \bullet  && &
}

\quad

Y:  \xymatrix{
\ar@{-}[d] & \ar@{-}[d] & & & & & & & & \\
\alpha \ar@{-}[d]\ar@{-}[r] & \beta \ar@{-}[d]\ar@{-}[r] & \bullet \ar@{-}[d]\ar@{-}[r] & \bullet \ar@{-}[d]\ar@{-}[r] & \bullet \ar@{-}[d] & & & & & \\
\bullet \ar@{-}[r] & \bullet \ar@{-}[r] & \bullet \ar@{-}[d]\ar@{-}[r] & \bullet \ar@{-}[d]\ar@{-}[r] & \bullet \ar@{-}[d] & & & & & \\
 & & \bullet \ar@{-}[r] & \bullet \ar@{-}[d]\ar@{-}[r] & \bullet \ar@{-}[d] & & & & & \\
 & & & \bullet \ar@{-}[d]\ar@{-}[r] & \bullet \ar@{-}[d]\ar@{-}[r]  & \bullet \ar@{-}[d]\ar@{-}[r] & \bullet \ar@{-}[d] & & & \\
 & & & \bullet \ar@{-}[r] & \bullet \ar@{-}[r]  & \bullet \ar@{-}[d]\ar@{-}[r] & \bullet \ar@{-}[d] & & & \\
 & & & & & \bullet \ar@{-}[d]\ar@{-}[r] & \bullet \ar@{-}[d]\ar@{-}[r]  & \bullet \ar@{-}[d] & & \\
 & & & & & \bullet \ar@{-}[d]\ar@{-}[r] & \bullet \ar@{-}[d]\ar@{-}[r] & \bullet \ar@{-}[d]\ar@{-}[r] & \bullet \ar@{-}[d]\ar@{-}[r]  & \bullet \ar@{-}[d] \\
 & & & & & \bullet \ar@{-}[r] & \bullet \ar@{-}[r] & \bullet \ar@{-}[r] & \gamma \ar@{-}[d]\ar@{-}[r]  & \delta \ar@{-}[d] \\
 & & & & & & & & & 
} 
\vskip 0.2cm
here  $\alpha$ and $\gamma$ are connected by an edge, and so are $\beta$ and $\delta$.

\quad

\quad

Z:  \xymatrix{
\ar@{-}[d] & \ar@{-}[d] & & & & & & & & \\
\alpha \ar@{-}[d]\ar@{-}[r] & \beta \ar@{-}[d]\ar@{-}[r] & \bullet \ar@{-}[d]\ar@{-}[r] & \bullet \ar@{-}[d] & & & & & & \\
\bullet \ar@{-}[r] & \bullet \ar@{-}[r] & \bullet \ar@{-}[d]\ar@{-}[r] & \bullet \ar@{-}[d] & & & & & & \\
 & & \bullet \ar@{-}[d]\ar@{-}[r] & \bullet \ar@{-}[d]\ar@{-}[r] & \bullet \ar@{-}[d]\ar@{-}[r]  & \bullet \ar@{-}[d] & & & & \\
 & & \bullet \ar@{-}[r] & \bullet \ar@{-}[r] & \bullet \ar@{-}[d]\ar@{-}[r]  & \bullet \ar@{-}[d] & & & & \\
 & & & & \bullet \ar@{-}[d]\ar@{-}[r] & \bullet \ar@{-}[d]\ar@{-}[r] & \bullet \ar@{-}[d]\ar@{-}[r]  & \bullet \ar@{-}[d] & & \\
 & & & & \bullet \ar@{-}[r] & \bullet \ar@{-}[r] & \bullet \ar@{-}[d]\ar@{-}[r]  & \bullet \ar@{-}[d] & & \\
 & & & & & & \bullet \ar@{-}[d]\ar@{-}[r] & \bullet \ar@{-}[d]\ar@{-}[r] & \bullet \ar@{-}[d]\ar@{-}[r]  & \bullet \ar@{-}[d] \\
 & & & & & & \bullet \ar@{-}[r] & \bullet \ar@{-}[r] & \gamma \ar@{-}[d]\ar@{-}[r]  & \delta \ar@{-}[d] \\
 & & & & & & & & &
} 
\vskip 0.2cm
here  $\alpha$ and $\gamma$ are connected by an edge, and so are $\beta$ and $\delta$.

\vskip 0.2cm

All these graphs can be embedded in the two or three dimensional integer lattice, even the graphs W, X, Y and Z as is shown below:
\vskip 1cm

W: \xy (24,52)*{} ; (24,-12)*{} **\dir{-} ,
(48,40)*{} ; (48,5)*{} **\dir{-} ,%
(48,3)*{} ; (48,-24)*{} **\dir{-} ,
(72,28)*{} ; (72,-36)*{} **\dir{-} ,
(-18,46)*{} ; (-18,-7)*{} **\dir{-} ,%
(-18,-9)*{} ; (-18,-18)*{} **\dir{-} ,
(54,10)*{} ; (54,-22)*{} **\dir{-} ,
(-36,28)*{} ; (-36,2)*{} **\dir{-} ,%
(-36,0)*{} ; (-36,-30)*{} **\dir{-} ,%
(-36,-32)*{} ; (-36,-36)*{} **\dir{-} ,
(-54,10)*{} ; (-54,-54)*{} **\dir{-} ,%
(-30,-2)*{} ; (-30,-34)*{} **\dir{-} ,
(24,-12)*{} ; (28,-14)*{} **\dir{-} ,
(30,-15)*{} ; (72,-36)*{} **\dir{-} ,
(-18,-18)*{} ; (30,-42)*{} **\dir{-} ,
(-36,-36)*{} ; (12,-60)*{} **\dir{-} ,
(24,20)*{} ; (50,7)*{} **\dir{-} ,%
(52,6)*{} ; (54,5)*{} **\dir{-} ,%
(56,4)*{} ; (72,-4)*{} **\dir{-} ,
(-54,-22)*{} ; (-30,-34)*{} **\dir{-} ,
(24,52)*{} ; (72,28)*{} **\dir{-} ,%
(-54,10)*{} ; (18,-26)*{} **\dir{-} ,
(-18,-18)*{} ; (-29,-29)*{} **\dir{-} ,
(-30,-30)*{} ; (-32,-32)*{} **\dir{-} ,
(-33,-33)*{} ; (-54,-54)*{} **\dir{-} ,
(24,-12)*{} ; (14,-22)*{} **\dir{-} ,%
(12,-24)*{} ; (-12,-48)*{} **\dir{-} ,
(48,-24)*{} ; (12,-60)*{} **\dir{-} ,
(-18,14)*{} ; (-32,0)*{} **\dir{-} ,
(-34,-2)*{} ; (-54,-22)*{} **\dir{-} ,
(72,-4)*{} ; (54,-22)*{} **\dir{-} ,
(-18,46)*{} ; (-54,10)*{} **\dir{-} ,
(72,28)*{} ; (18,-26)*{} **\dir{-} ,
\endxy
\quad
\vskip 1cm

X: \xy (24,52)*{} ; (24,26)*{} **\dir{-} ,%
(24,24)*{} ; (24,17)*{} **\dir{-} ,%
(24,15)*{} ; (24,-6)*{} **\dir{-} ,%
(24,-8)*{} ; (24,-12)*{} **\dir{-} ,
(48,40)*{} ; (48,-24)*{} **\dir{-} ,
(72,28)*{} ; (72,-36)*{} **\dir{-} ,
(-18,46)*{} ; (-18,20)*{} **\dir{-} ,%
(-18,18)*{} ; (-18,14)*{} **\dir{-} ,
(6,34)*{} ; (6,8)*{} **\dir{-} ,%
(6,6)*{} ; (6,2)*{} **\dir{-} ,
(30,22)*{} ; (30,-10)*{} **\dir{-} ,
(-12,-16)*{} ; (-12,-42)*{} **\dir{-} ,%
(-12,-44)*{} ; (-12,-48)*{} **\dir{-} ,
(12,-28)*{} ; (12,-60)*{} **\dir{-} ,
(36,-40)*{} ; (36,-72)*{} **\dir{-} ,
(-54,10)*{} ; (-54,-54)*{} **\dir{-} ,%
(-30,-2)*{} ; (-30,-66)*{} **\dir{-} ,
(-6,-14)*{} ; (-6,-78)*{} **\dir{-} ,
(24,-12)*{} ; (72,-36)*{} **\dir{-} ,
(6,-30)*{} ; (8,-31)*{} **\dir{-} ,%
(10,-32)*{} ; (12,-33)*{} **\dir{-} ,%
(14,-34)*{} ; (34,-44)*{} **\dir{-} ,%
(38,-46)*{} ; (54,-54)*{} **\dir{-} ,
(-12,-48)*{} ; (-8,-50)*{} **\dir{-} ,%
(-4,-52)*{} ; (36,-72)*{} **\dir{-} ,
(-54,-54)*{} ; (-6,-78)*{} **\dir{-} ,
(24,20)*{} ; (26,19)*{} **\dir{-} ,
(28,18)*{} ; (30,17)*{} **\dir{-} ,
(32,16)*{} ; (72,-4)*{} **\dir{-} ,
(-18,14)*{} ; (30,-10)*{} **\dir{-} ,
(-12,-16)*{} ; (-8,-18)*{} **\dir{-} ,%
(-4,-20)*{} ; (36,-40)*{} **\dir{-} ,
(-54,-22)*{} ; (-6,-46)*{} **\dir{-} ,
(24,52)*{} ; (72,28)*{} **\dir{-} ,%
(-18,46)*{} ; (30,22)*{} **\dir{-} ,
(-36,28)*{} ; (12,4)*{} **\dir{-} ,
(-54,10)*{} ; (-6,-14)*{} **\dir{-} ,
(6,-30)*{} ; (-5,-41)*{} **\dir{-} ,
(-6,-42)*{} ; (-8,-44)*{} **\dir{-} ,
(-10,-46)*{} ; (-30,-66)*{} **\dir{-} ,
(30,-42)*{} ; (-6,-78)*{} **\dir{-} ,
(72,-36)*{} ; (36,-72)*{} **\dir{-} ,
(24,20)*{} ; (10,6)*{} **\dir{-} ,
(8,4)*{} ; (6,2)*{} **\dir{-} ,
(-12,-16)*{} ; (-30,-34)*{} **\dir{-} ,
(48,8)*{} ; (30,-10)*{} **\dir{-} ,
(12,-28)*{} ; (-6,-46)*{} **\dir{-} ,
(-18,46)*{} ; (-54,10)*{} **\dir{-} ,
(24,52)*{} ; (-12,16)*{} **\dir{-} ,
(48,40)*{} ; (12,4)*{} **\dir{-} ,
\endxy
\quad
\vskip 1cm

Y: \xy (0,64)*{} ; (0,29)*{} **\dir{-} ,%
(0,27)*{} ; (0,11)*{} **\dir{-} ,%
(0,9)*{} ; (0,0)*{} **\dir{-} ,
(24,52)*{} ; (24,26)*{} **\dir{-} ,%
(24,24)*{} ; (24,17)*{} **\dir{-} ,%
(24,15)*{} ; (24,-1)*{} **\dir{-} ,%
(24,-3)*{} ; (24,-12)*{} **\dir{-} ,
(48,40)*{} ; (48,14)*{} **\dir{-} ,%
(48,12)*{} ; (48,5)*{} **\dir{-} ,%
(48,3)*{} ; (48,-24)*{} **\dir{-} ,
(-18,14)*{} ; (-18,-7)*{} **\dir{-} ,%
(-18,-9)*{} ; (-18,-18)*{} **\dir{-} ,
(54,10)*{} ; (54,-22)*{} **\dir{-} ,
(-36,-4)*{} ; (-36,-36)*{} **\dir{-} ,%
(36,-8)*{} ; (36,-40)*{} **\dir{-} ,
(-30,-2)*{} ; (-30,-66)*{} **\dir{-} ,%
(-6,-14)*{} ; (-6,-78)*{} **\dir{-} ,
(18,-26)*{} ; (18,-90)*{} **\dir{-} ,
(0,0)*{} ; (28,-14)*{} **\dir{-} ,
(30,-15)*{} ; (34,-17)*{} **\dir{-} ,
(38,-19)*{} ; (48,-24)*{} **\dir{-} ,
(-18,-18)*{} ; (-8,-23)*{} **\dir{-} ,%
(-4,-25)*{} ; (16,-35)*{} **\dir{-} ,%
(20,-37)*{} ; (30,-42)*{} **\dir{-} ,
(-36,-36)*{} ; (-32,-38)*{} **\dir{-} ,%
(-28,-40)*{} ; (-8,-50)*{} **\dir{-} ,%
(-4,-52)*{} ; (12,-60)*{} **\dir{-} ,
(-30,-66)*{} ; (18,-90)*{} **\dir{-} ,
(0,32)*{} ; (2,31)*{} **\dir{-} ,%
(4,30)*{} ; (26,19)*{} **\dir{-} ,%
(28,18)*{} ; (48,8)*{} **\dir{-} ,
(-30,-34)*{} ; (18,-58)*{} **\dir{-} ,
(0,64)*{} ; (48,40)*{} **\dir{-} ,%
(6,34)*{} ; (54,10)*{} **\dir{-} ,
(-12,16)*{} ; (36,-8)*{} **\dir{-} ,
(-30,-2)*{} ; (18,-26)*{} **\dir{-} ,
(0,0)*{} ; (-10,-10)*{} **\dir{-} ,
(-12,-12)*{} ; (-29,-29)*{} **\dir{-} ,
(-31,-31)*{} ; (-36,-36)*{} **\dir{-} ,
(6,-30)*{} ; (-5,-41)*{} **\dir{-} ,%
(-7,-43)*{} ; (-9,-45)*{} **\dir{-} ,%
(-10,-46)*{} ; (-30,-66)*{} **\dir{-} ,
(30,-42)*{} ; (19,-53)*{} **\dir{-} ,
(17,-55)*{} ; (15,-57)*{} **\dir{-} ,
(14,-58)*{} ; (-6,-78)*{} **\dir{-} ,
(0,32)*{} ; (-36,-4)*{} **\dir{-} ,
(54,-22)*{} ; (18,-58)*{} **\dir{-} ,
(24,52)*{} ; (-12,16)*{} **\dir{-} ,
(48,40)*{} ; (12,4)*{} **\dir{-} ,
(54,10)*{} ; (18,-26)*{} **\dir{-} ,
\endxy
\quad
\vskip 1cm

Z: \xy (0,64)*{} ; (0,38)*{} **\dir{-} ,%
(0,36)*{} ; (0,29)*{} **\dir{-} ,%
(0,27)*{} ; (0,11)*{} **\dir{-} ,%
(0,9)*{} ; (0,0)*{} **\dir{-} ,
(24,52)*{} ; (24,26)*{} **\dir{-} ,
(24,24)*{} ; (24,17)*{} **\dir{-} ,
(24,15)*{} ; (24,-12)*{} **\dir{-} ,
(-18,46)*{} ; (-18,20)*{} **\dir{-} ,%
(-18,18)*{} ; (-18,14)*{} **\dir{-} ,
(-36,28)*{} ; (-36,-4)*{} **\dir{-} ,%
(36,-8)*{} ; (36,-72)*{} **\dir{-} ,
(-30,-2)*{} ; (-30,-34)*{} **\dir{-} ,%
(-6,-14)*{} ; (-6,-46)*{} **\dir{-} ,
(18,-26)*{} ; (18,-90)*{} **\dir{-} ,
(0,0)*{} ; (4,-2)*{} **\dir{-} ,
(6,-3)*{} ; (28,-14)*{} **\dir{-} ,
(30,-15)*{} ; (34,-17)*{} **\dir{-} ,
(38,-19)*{} ; (72,-36)*{} **\dir{-} ,
(-18,-18)*{} ; (-8,-23)*{} **\dir{-} ,
(-4,-25)*{} ; (16,-35)*{} **\dir{-} ,
(20,-37)*{} ; (32,-43)*{} **\dir{-} ,
(34,-44)*{} ; (36,-45)*{} **\dir{-} ,
(38,-46)*{} ; (54,-54)*{} **\dir{-} ,
(12,-60)*{} ; (16,-62)*{} **\dir{-} ,
(20,-64)*{} ; (36,-72)*{} **\dir{-} ,
(-6,-78)*{} ; (18,-90)*{} **\dir{-} ,
(0,32)*{} ; (2,31)*{} **\dir{-} ,%
(4,30)*{} ; (24,20)*{} **\dir{-} ,
(-30,-34)*{} ; (18,-58)*{} **\dir{-} ,
(0,64)*{} ; (24,52)*{} **\dir{-} ,%
(-18,46)*{} ; (30,22)*{} **\dir{-} ,
(-36,28)*{} ; (12,4)*{} **\dir{-} ,
(-30,-2)*{} ; (18,-26)*{} **\dir{-} ,
(0,0)*{} ; (-10,-10)*{} **\dir{-} ,
(-12,-12)*{} ; (-18,-18)*{} **\dir{-} ,
(24,-12)*{} ; (14,-22)*{} **\dir{-} , %
(12,-24)*{} ; (6,-30)*{} **\dir{-} , 
(48,-24)*{} ; (37,-35)*{} **\dir{-} , 
(35,-37)*{} ; (19,-53)*{} **\dir{-} , 
(17,-55)*{} ; (15,-57)*{} **\dir{-} , 
(14,-58)*{} ; (-6,-78)*{} **\dir{-} , 
(72,-36)*{} ; (18,-90)*{} **\dir{-} ,
(0,32)*{} ; (-14,18)*{} **\dir{-} ,
(-16,16)*{} ; (-36,-4)*{} **\dir{-} ,
(36,-40)*{} ; (18,-58)*{} **\dir{-} ,
(0,64)*{} ; (-36,28)*{} **\dir{-} ,
(6,34)*{} ; (-30,-2)*{} **\dir{-} ,
(30,22)*{} ; (-6,-14)*{} **\dir{-} ,
(36,-8)*{} ; (18,-26)*{} **\dir{-} ,
\endxy
\quad

\vskip 1cm


\begin{thebibliography}{1}
\bibitem[1]{bb}  G. Bowlin, M. G.  Brin, {\it Coloring planar graphs via colored paths in the associahedra}, 
Internat. J. Algebra Comput. {\bf 23} (2013), no. 6, 1337-1418.
\bibitem[2]{c}  R. P. Carpentier, {\it On signed diagonal flip sequences}, European J. Combin. {\bf 32} (2011), no. 3,  472-477.
\bibitem[3]{d} D. {\u{Z}}. Djokovi{\'{c}}, {\it Distance-preserving subgraphs of hypercubes}, J. Combinatorial Theory Ser. B 14 (1973) 263-267. 
\bibitem[4]{e}  S. Eliahou, {\it Signed diagonal flips and the four color theorem}, European J. Combin. {\bf 20} (1999), no. 7,  641-647.
\bibitem[5]{el}  S. Eliahou, C. Lecouvey, {\it Signed permutations and the four color theorem}, Expo. Math. {\bf 27} (2009), no. 4, 313-340.
\bibitem[6]{gp}  S. Gravier, C. Payan, {\it Flips sign\'es et triangulations d'un polygone}, European J. Combin. {\bf 23} (2002), no. 7,  817-821.
\bibitem[7]{o} S. Ovchinnikov, {\it Partial cubes: structures, characterizations, and constructions}, Discrete Math. {\bf 308} (2008), no. 23,   5597-5621.






\end{thebibliography}
\end{document}